\documentclass[11pt]{amsart}
\usepackage{amsfonts,amssymb,amscd,amsmath,enumerate,verbatim,calc, hyperref, thumbpdf}

\theoremstyle{plain}
\newtheorem{Theorem}{Theorem}[section]

\newtheorem{Corollary}[Theorem]{Corollary} 
\newtheorem{Proposition}[Theorem]{Proposition}

\theoremstyle{definition}
\newtheorem{Definition}[Theorem]{{Definition}}
\newtheorem{Example}[Theorem]{Example}

\newtheorem*{chunk*}{}
\newtheorem*{ack}{Acknowledgement}
\numberwithin{equation}{section}

\theoremstyle{remark}
\newtheorem{Remark}[Theorem]{Remark}

\newcommand{\fm}{{\mathfrak m}}

\newcommand{\fp}{{\mathfrak p}}

\usepackage{graphicx}

\newcommand{\showtwo}[3]
{%
    \begin{center}
        \includegraphics[width=#3\linewidth]{#1}%
        \hfill
        \includegraphics[width=#3\linewidth]{#2}%
    \end{center}
}

\begin{document}
\subjclass{Primary 13A15. Secondary 13A18}
\title[Asymptotic growth of powers of ideals]
{Asymptotic growth of powers of ideals}

\author{C\u{a}t\u{a}lin Ciuperc\u{a}}
\address{Department of Mathematics, North Dakota State University, Fargo, ND 58105}
\email{catalin.ciuperca@ndsu.edu}

\author{Florian Enescu}
\address{Department of Mathematics and Statistics, Georgia State University, Atlanta, GA 30303}
\email{fenescu@mathstat.gsu.edu}
\thanks{The second author gratefully acknowledges partial financial support from the National Science Foundation, CCF-0515010 and Georgia State 
University, Research Initiation Grant.}

\author{Sandra Spiroff}
\address{Department of Mathematics, Seattle University, Seattle, WA  98122} 
\email{spiroffs@seattleu.edu}

%\subjclass{13}
%\subjclass[2000]{Primary 13-XX}

\begin{abstract} 
Let $A$ be a locally analytically unramified local ring and  $J_1,\ldots, J_k, I$  ideals such that $J_i \subseteq \sqrt{I}$ for all $i$, the ideal $I$ is not nilpotent, and $\bigcap_k I^k=(0)$. Let  $C=C(J_1,\ldots,J_k;I) \subseteq \mathbb{R}^{k+1}$ be the cone generated by $\{ (m_1,\ldots,m_k,n) \in \mathbb{N}^{k+1} \mid J_1^{m_1}\ldots J_k^{m_k} \subseteq I^n \}$.  We prove that the topological closure of $C$ is a rational  polyhedral cone. This generalizes results by Samuel, Nagata, and Rees.
\end{abstract}
\maketitle
%\bigskip

\section*{Introduction} In this note we continue the  study of the asymptotic properties of powers of ideals initiated by Samuel in \cite{Samuel1}. Let $A$ be a commutative noetherian ring with identity and $I, J$ ideals in $A$ with $J \subseteq \sqrt{I}$. Also, assume that the ideal $I$ is not nilpotent and $\bigcap_k I^k=(0)$. Then for each positive integer $m$ one can define $v_I(J,m)$ to be the largest integer $n$ such that $J^m \subseteq I^n$. Similarly, $w_J(I,n)$ is defined to be the smallest integer $m$ such that $J^m \subseteq I^n$. Under the above assumptions, Samuel  proved that the sequences  $\{v_I(J,m)/m\}_m$ and $\{w_J(I,n)/n\}_n$ have limits $l_I(J)$ and $L_J(I)$, respectively,   and  $l_I(J) L_J(I)=1$ \cite[Theorem 1]{Samuel1}. It is also observed that these limits are actually the supremum  and infimum of the respective sequences. One of the questions raised in Samuel's paper is whether $l_I(J)$ is always rational. This has been positively answered  by Nagata \cite{Nagata1} and Rees \cite{Rees1}. The approach used by Rees is described in the next  section of this paper.  

We consider the following generalization of the problem described
above. Let  $J_1,\ldots, J_k, I$ be ideals in a locally analytically
unramified ring $A$ such that $J_i \subseteq \sqrt{I}$ for all $i$,
$I$ is not nilpotent, and $\bigcap_k I^k=(0)$, and let
$C=C(J_1,\ldots,J_k;I) \subseteq \mathbb{R}^{k+1}$ be the cone
generated by $\{ (m_1,\ldots,m_k,n) \in \mathbb{N}^{k+1} \mid
J_1^{m_1}\ldots J_k^{m_k} \subseteq I^n \}.$ We prove that the
topological closure of $C$ is a rational polyhedral cone; i.e., a
polyhedral cone bounded by hyperplanes whose equations have rational
coefficients. Note that the case $k=1$  follows from the results
proved by Samuel, Nagata, and Rees; the cone $C$ is the intersection of
the half-planes given by $n \geq 0$ and $ n \leq l_I(J) m_1$. In
Section 3 we look at the periodicity of the rate of change of the sequence $\{v_I(J,m)\}_m$, more precisely, the 
periodicity of the sequence $\{v_I(J,m+1) - v_I(J,m) \}_m$. The last part of the paper describes a method  of computing the limits studied by Samuel in the case of monomial ideals.

\section{The Rees valuations of an ideal}

In this section we give a brief description of the Rees valuations associated to an ideal. 

For a noetherian ring $A$ which is not necessarily an integral domain, a discrete  valuation on $A$ is defined as follows.

\begin{Definition} Let $A$ be a noetherian ring. We say that $v: A \to \mathbb{Z} \cup \{\infty \}$ is a discrete valuation on $A$ if $\{x \in A \mid v(x)=\infty\}$ is a prime ideal $P$, $v$ factors through $A \to A/P \to \mathbb{Z} \cup \{ \infty \}$, and the induced function on $A/P$ is a rank one discrete valuation on $A/P$. If $I$ is an ideal in $A$, then we denote $v(I):=\min \{v(x) \mid x \in I \}$. 
\end{Definition}

If $R$ is a noetherian ring, we  denote by $\overline{R}$ the integral closure of $R$ in its total quotient ring $Q(R)$. 

\begin{Definition}
Let $I$ be an ideal in a noetherian ring $A$. An element $x \in A$ is
said to be integral over $I$ if $x$ satisfies an equation $x^n + a_1
x^{n-1}+\ldots+a_n=0$ with $a_i \in I^i$. The set of all elements in
$A$ that are integral over $I$ is an ideal $\overline{I}$, and the
ideal $I$ is called integrally closed if $I=\overline{I}$.  If all the
powers $I^n$ are integrally closed, then $I$ is said to be normal.
\end{Definition}

Given an ideal   $I$ in a noetherian ring $A$,  for each $x \in A$ let $v_I (x)=\sup \{ n \in \mathbb{N} \mid x \in I^n \}$. Rees \cite{Rees1} proved that for each $x \in A$ one can define  $$\overline{v}_I(x)=\lim_{k \to \infty} \frac{v_I(x^k)}{k},$$  and  for each integer $n$ one has $\overline{v}_I(x) \geq n$ if and only if $x \in \overline{I^n}$. Moreover, there exist discrete valuations $v_1,\ldots, v_h$ on $A$ in the sense defined above, and positive integers $e_1, \ldots, e_h$ such that, for each $x \in A$, 

\begin{equation}\label{num1}
\overline{v}_I(x)=\min \Big\{ \frac{v_i(x)}{e_i} \mid i=1,\ldots, h \Big\}.
\end{equation}

We briefly describe a construction of the Rees valuations $v_1,\ldots, v_h$. Let  $\fp_1, \ldots, \fp_g$ be the minimal prime ideals $\fp$ in $A$ such that $\fp+I \neq A$, and  let $\mathcal{R}_i(I)$ be the Rees ring $(A/\fp_i)[It, t^{-1}]$. Denote by  $W_{i1},\ldots, W_{ih_i}$  the rank one discrete valuation rings obtained by  localizing the rings $\overline{\mathcal{R}_i(I)}$
 at the minimal primes over $t^{-1}\overline{\mathcal{R}_i(I)}$,  let $w_{ij}$ ($i=1,\ldots, g$, $1 \leq j \leq h_i$)  be the corresponding discrete valuations, and let $V_{ij}=W_{ij} \cap Q(A/\fp_i)$ ($i=1,\ldots,g$).  
Then define $v_{ij}(x):=w_{ij}(x+\fp_i)$ and  $e_{ij}:=w_{ij}(t^{-1})(=v_{ij}(I))$ for all $i$, and for simplicity, renumber them as $e_1,\ldots, e_h$ and $v_1,\ldots, v_h$, respectively.

 Rees \cite{Rees1} proved that $v_1,\ldots,v_h$  are valuations satisfying (\ref{num1}). We refer the reader to the original article \cite{Rees1} for more details on this construction. 

\begin{Remark} With the notation established above, for every positive integer $n$ we have $$\overline{I^n}=\bigcap_{i=1}^{h} I^n V_i \cap R.$$ 
\end{Remark}

In particular, we have the following.
\begin{Remark}\label{rem2}If $K,L$ are ideals in $A$, $v_1,\ldots,v_h$ are the Rees valuations of $L$, and $v_i(K) \geq v_i(L)$ for all $i=1,\ldots,h$, then $\overline{K} \subseteq \overline{L}$.
\end{Remark}

The rationality of $l_I(J)$ can now be obtained as  consequence of the results of Rees. Indeed, by \cite[Theorem 2]{Samuel1}, if $J=(a_1,\ldots a_s)$, then $l_I(J)=\min \{l_I(a_i) \mid i=1,\ldots s \}$, and for each $i$ we have  $l_I(a_i)=\overline{v}_I(a_i)$, which is rational.

Finally, recall the following definition.
\begin{Definition} A local noetherian ring $(A, \fm)$ is analytically unramified if its $\fm$-adic completion $\hat A$ is reduced.
\end{Definition}

Rees \cite{Rees2} proved that for every ideal $I$ in an analytically unramified ring there exists an integer $k$ such that for all $n \geq 0$, $\overline{I^{n+k}} \subseteq I^n$.
\section{The cone structure}

Throughout this section $A$ is a locally analytically unramified ring and $I$ and $\underline{J}=J_1,\ldots, J_k$ are ideals in $A$ such that $J_i \subseteq \sqrt{I}$ for all $i$. Let $C=C(J_1,\ldots,J_k;I)\subseteq \mathbb{R}^{k+1}$ denote the cone generated by $\{ (m_1,\ldots,m_k,n) \in \mathbb{N}^{k+1} \mid J_1^{m_1} \ldots J_k^{m_k} \subseteq I^n \}.$ Also, for  $(m_1,\ldots,m_k)\in \mathbb{N}^{k}$, let  $v_I(\underline{J},m_1,\ldots,m_k)$ denote the largest nonnegative integer $n$ such that $ J_1^{m_1} \ldots J_k^{m_k} \subseteq I^n$.

For each Rees valuation $v_j$ of $I$, denote
$\alpha_{ij}=v_j(J_i)/e_j$ for all $i,j$, where $e_j = v_j(I)$.  Then we consider

$$D_j=\{(m_1,\ldots,m_k) \in {\mathbb{R}}_{\geq 0}^k \mid  \sum_{s=1}^{k} m_s \alpha_{sj} \leq  \sum_{s=1}^{k} m_s \alpha_{sl} \text{ for all } l \neq j \},$$ 
and we say that a Rees valuation $v_j$ is relevant if $D_j \neq \{0\}$. After a renumbering, assume that $v_1,v_2, \ldots, v_r$ ($r \leq h$) are the relevant Rees valuations. 

Note that each $D_j$ is an intersection of half-spaces (hence a polyhedral cone),   $\bigcup_{j=1}^{r} D_j = \mathbb{R}_{ \geq 0}^{k}$,  and two cones  $D_i, D_j$ ($ i \neq j$) either intersect along one common face or have only the origin in common. 
Let $$E_j= \{ (m_1,\ldots,m_k,n) \in {\mathbb{R}}_{+}^{k+1} \mid (m_1,\ldots,m_k) \in D_j \text{ and } n <  \sum_{s=1}^{k} m_s \alpha_{sj} \}$$ and $$\overline{E}_j=\{ (m_1,\ldots,m_k,n) \in {\mathbb{R}}_{+}^{k+1} \mid (m_1,\ldots,m_k) \in D_j \text{ and } n \leq  \sum_{s=1}^{k} m_s \alpha_{sj} \}.$$

\begin{Theorem} Let $A$ be a locally  analytically unramified ring. 
  Then for
  each $j=1,\ldots,r$ we have $$ E_j \cap \mathbb{Q}^{k+1} \subseteq C \cap (D_j
  \times \mathbb{R}_{\geq 0}) \subseteq \overline{E}_j.$$ 
\end{Theorem}

\begin{proof}  Let $(m_1,\ldots,m_k,n) \in  C \cap (D_j \times
  \mathbb{R}_{\geq 0})$. Then there exists $t \in \mathbb{R}$ such
  that $tm_1,\ldots, tm_k$ are positive integers and
  $$J_1^{tm_1}\ldots J_k^{tm_k} \subseteq I^{tn}.$$  Hence, for each  Rees valuation $v_j$ of $I$ we obtain $$ tm_1v_j(J_1)+\cdots +tm_k v_j(J_k) \geq tnv_j(I),$$ or  equivalently, $$ n \leq \sum_{s=1}^{k} m_s \alpha_{sj}.$$

For the other inclusion, first observe that it is enough to prove that  $E_j \cap \mathbb{Z}^{k+1} \subseteq C \cap (D_j \times \mathbb{R}_{\geq 0})$. Indeed, if $E_j \cap \mathbb{Z}^{k+1}\subseteq C \cap (D_j \times \mathbb{R}_{\geq 0}) $, then for each $\alpha \in E_j \cap \mathbb{Q}^{k+1}$ there exists a positive integer $L$ such that $\alpha L \in E_j \cap \mathbb{Z}^{k+1} \subseteq C \cap (D_j \times \mathbb{R}_{\geq 0})$. This implies that $\alpha \in (1/L)\big( C \cap (D_j \times \mathbb{R}_{\geq 0}) \big)=C \cap (D_j \times \mathbb{R}_{\geq 0})$

Let $(m_1,\ldots,m_k,n) \in  E_j \cap \mathbb{Z}^{k+1}$.  Set $\alpha=\sum_{s=1}^{k} m_s \alpha_{sj}$. Since the ring $A$ is  analytically unramified, there exists an integer $N$ such that $\overline{I^t} \subseteq I^{t-N}$ for all $t$. (The convention is that $I^n=A$ for $n \leq 0$.) Let $g$ be the integer part of $\alpha$. For any Rees valuation $v_i$ of $A$ we then get $$v_i(I^g)=ge_i \leq \alpha e_i \leq (\sum_{s=1}^{k} m_s \alpha_{si}) e_i = v_i(J_1^{m_1}\ldots J_k^{m_k}),$$ and hence, by Remark \ref{rem2}, $$J_1^{m_1}\ldots J_k^{m_k} \subseteq \overline{I^g} \subseteq I^{g-N}.$$ This implies that 
\begin{equation}\label{eq2}
v_I(\underline{J},m_1,\ldots,m_k) \geq g-N > \alpha - 1 -N. 
\end{equation}

Since  $ n < \alpha$, we can find $\delta >0$ such that $n < \alpha -\delta$. Choose $l$ such that $l\delta > N+1$ and $l m_1, \ldots, l m_k, l n$ are integers. By (\ref{eq2}),  we obtain $v_I(\underline{J}, l m_1,\ldots,l m_k) > l\alpha -N -1$, and by the choice of $l$, we also have  $n l < l\alpha -N -1$. Then  $nl < v_I(\underline{J}, l m_1,\ldots,l m_k)$, which  implies that $J_1^{l m_1}\ldots J_k^{l m_k} \subseteq I^{l n}$; i.e., $(m_1,\ldots,m_k,n) \in C$.
\end{proof}

\begin{Corollary} The topological closure of $C$ is a rational polyhedral cone.
\end{Corollary}

\begin{proof} From the previous theorem it follows that the
  topological closure of $C \cap (D_j \times \mathbb{R}_{\geq 0})$ is
  $\overline{E}_j$, and hence the topological closure of $C$ is the polyhedral cone bounded by the hyperplanes $n=\sum_{s=1}^{k}m_s\alpha_{sj}$ ($j=1,\ldots, r$) and the coordinate hyperplanes.
\end{proof}

A detailed example of Corollary 2.2 is given below in Example 2.5.

\begin{Corollary} Let $a_1,a_2,\ldots, a_k$ be real numbers. The limit 
\begin{equation}\label{eq-c}
\lim_{m_1,\ldots,m_k \to \infty} \frac{v_I(\underline{J},m_1,\ldots,m_k)}{a_1m_1+\ldots+a_k m_k}
\end{equation}
 exists if and only if there exists a rational number $l$ such that  $la_s=\alpha_{s1}=\alpha_{s2}=\ldots=\alpha_{sr}$ for all  $s=1,\ldots,k$. In this case the limit is equal to $l$.
\end{Corollary} 

\begin{proof} Since the polyhedral cones $D_j$ form a partition of $\mathbb{R}_{\geq 0}^k$,  the limit (\ref{eq-c}) exists and is equal to $l$ if and only if for each $j$ we have
\begin{equation}\label{eq3}
\lim_{\substack{m_1,\ldots,m_k \to \infty \\ (m_1,\ldots,m_k) \in D_j}} \frac{v_I(\underline{J},m_1,\ldots,m_k)}{a_1m_1+\ldots+a_k m_k}=l.
\end{equation}
On the other hand, (\ref{eq3}) holds if  and only if 
$l a_s=\alpha_{sj}$ for all $s=1,\ldots, k$. Indeed, this limit exists  and is equal to $l$ if and only if 
over $D_j$  the topological closure of $C$ is bounded by the hyperplane $n=la_1m_1+ \ldots+la_k m_k$, which therefore should coincide with the hyperplane $n=\sum_{s=1}^{k}m_s\alpha_{sj}$. 

In conclusion, the limit (\ref{eq-c}) exists and is equal to $l$ if and only if all the hyperplanes $n=\sum_{s=1}^{k}m_s\alpha_{sj}$ ($j=1,\ldots,r$) coincide with $n=la_1m_1+ \ldots+la_k m_k$, or equivalently, $la_s=\alpha_{s1}=\alpha_{s2}=\ldots=\alpha_{sr}$ for all  $s=1,\ldots,k$.
\end{proof}

\begin{Corollary} Assume that the ideal $I$ has only one Rees valuation. Then 
 the limit  $$\lim_{m_1,\ldots,m_k \to \infty} \frac{v_I(\underline{J},m_1,\ldots,m_k)}{a_1m_1+\ldots+a_k m_k}$$ exists if and only if $ l_I(J_1)/ a_1=\ldots= l_I(J_k)/ a_k$.
\end{Corollary}

\begin{proof} This is a particular case of the previous Corollary.
\end{proof}

\begin{Example}  Let $A = \displaystyle{\mathbb R[[X,Y,Z]]/(XY^2 -
  Z^9)}$ and $I = (x,y,z)R$ as in \cite[Example 3.1]{HS}.  Then $\mathcal R(I) = A[It,t^{-1}]$,  $\displaystyle{\mathcal R(I)/t^{-1} \mathcal R(I)}$
  $\displaystyle{\cong Q[xt,yt,zt]/(xt)(yt)}$, and there are two Rees
  valuations $v_1$ and $v_2$, corresponding to the minimal primes
  $\mathfrak p_1 = (xt,t^{-1})$ and $\mathfrak p_2 = (yt,t^{-1})$, over
  $t^{-1} \mathcal R(I)$.  As shown in  \cite[Example 3.1]{HS}, we have $v_1(x)=7,
  v_1(y)=v_1(z)=1$ and $v_2(x)=v_2(z)=1, v_2(y)=4$.  Thus $v_1(I) =
  \min\{v_1(x), v_1(y), v_1(z)\} = 1$.  Likewise $v_2(I) = 1$.  Set
  $J_1= (x,z^2)$ and $J_2 = (y^2,z^3)$.  Then $v_1(J_1) =2, v_2(J_1)
  = 1$, and $v_1(J_2) = 2, v_2(J_2) = 3$. 
  Therefore, $E_1 = \{(m_1,m_2,n) | n \leq 2m_1 + 2m_2 \}$ and $E_1 =
  \{(m_1,m_2,n) | n \leq m_1 + 3m_2 \}$.  The boundary planes of
  $E_1$ and $E_2$ in $\mathbb R^3$ are $z = 2x + 2y$ and $z=x+3y$, respectively.  Thus, 
  according to the results of Corollary 2.2, the
  topological closure of the cone generated by $\{(m_1,m_2,n) |
  J_1^{m_1}J_2^{m_2} \subseteq I^n \}$ is as pictured below. 
\end{Example}

\begin{figure}[tbh]
    \showtwo{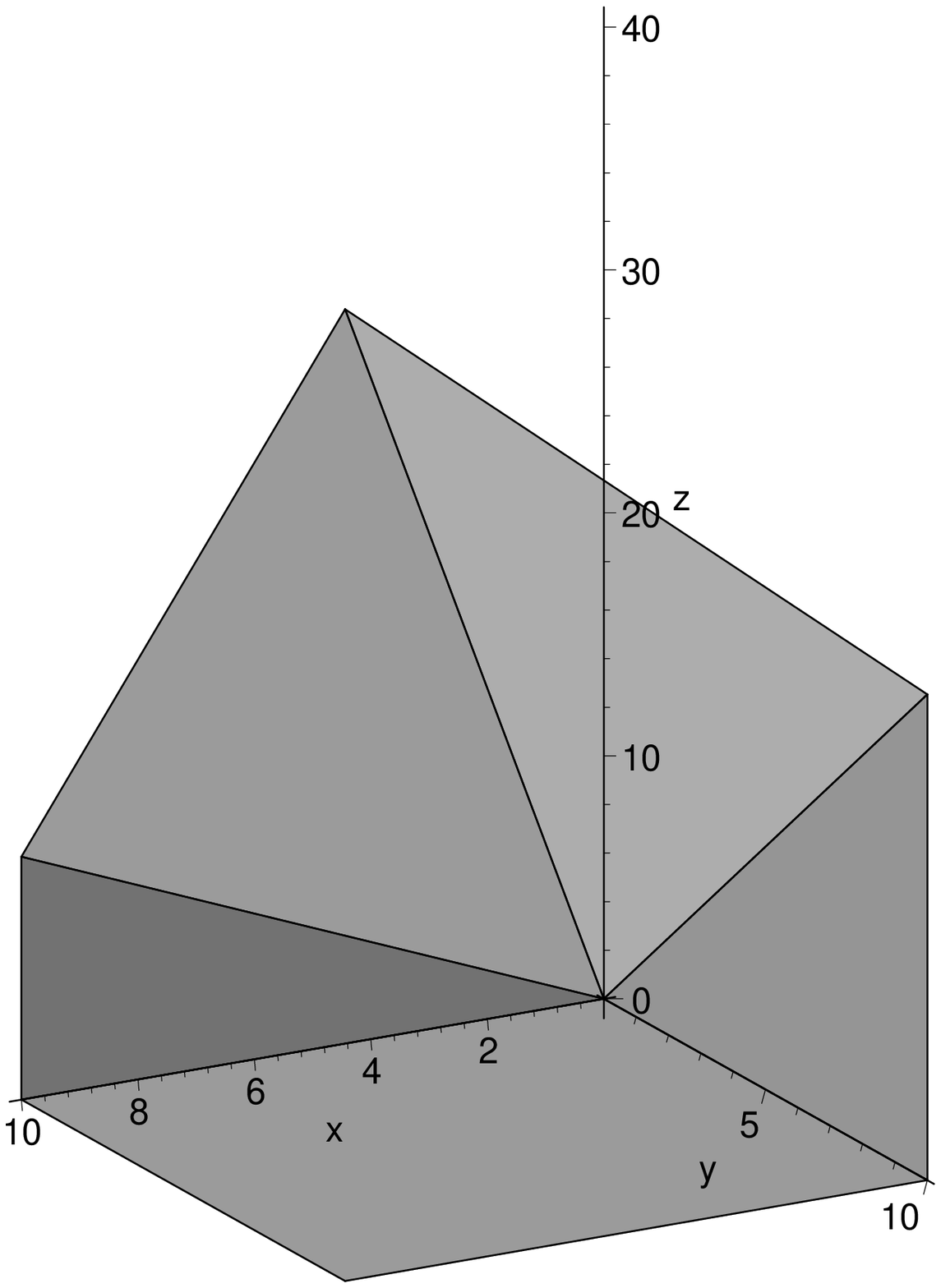}{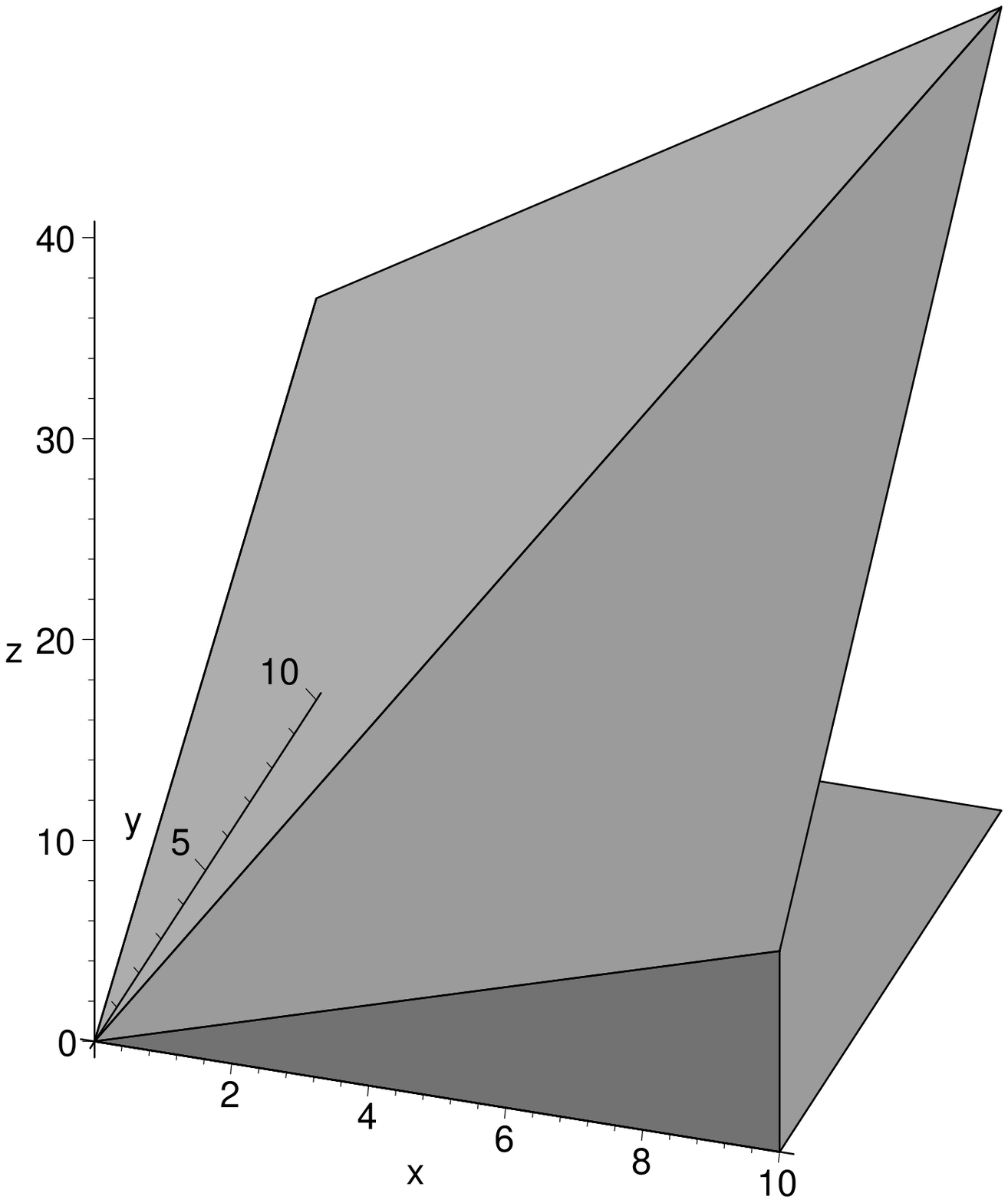}{0.4}
    \caption{View from the front and rotated $90^{\circ}$
    ctr-clockwise around the $z$-axis.}
\end{figure}

\begin{Example}  Let $A = \displaystyle{k[[X,Y]]}$, with $k$ a
    field, and $I = (x^3, x^2y, y^2)$.  As shown in \cite{Sally},
    $I$ has only one associated Rees valuation.  Let $J_1
  = (x^3y^7), J_2 = (x^4y^6)$, and $J_3 =(x^5y^2)$.  Using the methods
  in Section 4, we can compute $\displaystyle{l_I(J_1)= 9/2, l_I(J_2)=
  13/3}$, and $\displaystyle{l_I(J_3)= 8/3}$.  Then by
  Corollary 2.4, 
  the limit  $$\lim_{m_1,m_2,m_3 \to \infty} 
\frac{v_I(J_1,J_2,J_3,m_1,m_2,m_3)}{27m_1+26m_2+16 m_3}$$ exists and equals $\displaystyle{\frac{1}{6}}$ since
$\displaystyle{\frac{l_I(J_1)}{27} =
  \frac{l_I(J_2)}{26} =
  \frac{l_I(J_3)}{16}} = \frac{1}{6}$.
\end{Example}

\section{Periodic Increase}

In this section we take a closer look at  the sequence $\{ v_I(J,m)\}_m$. To simplify the notation we will simply write $v(m)$ instead of $v_I(J,m)$. 

We address the question of whether this sequence increases eventually in a periodic way; that is, whether or not there exists a positive integer $t$ such that 
$v(m+t)-v(m+t-1)=v(m)-v(m-1)$ for $ m \gg 0$, or equivalently, $v(m+t)-v(m)=constant$, for $m \gg 0$. Our work is partly motivated by \cite[Theorem 8]{Nagata1}, where Nagata proves that the deviation $v(m)-l_I(J)m$ is bounded. In particular, this implies that  there exists a positive constant $C$ such that $0 \leq v(m+t)-v(m)-v(t) < C$ for all $m,t$.

We begin by defining noetherian filtrations.

\begin{Definition} A family of ideals $\mathcal{F}=\{F_m\}_{m \geq 0}$ in a noetherian ring  $A$ is called a filtration if $F_0=A$, $F_{m+1} \subseteq F_m$, and $F_m F_n \subseteq F_{m+n}$ for all $ m,n \geq 0$. We say that the  filtration $\{F_m\}_{m \geq 0}$ is noetherian if the associated graded ring $\oplus_{m \geq 0} F_m$ is noetherian. Equivalently, the filtration  $\mathcal{F}$ is noetherian   if and only if there exists $t$ such that $F_{m+t}=F_m F_t$ for all $m \geq t$ (\cite[4.5.12]{BH}).
\end{Definition}

\begin{Proposition}\label{prop-per} Let $I,J$ be ideals in a noetherian local ring $A$ such that  $J \subseteq \sqrt{I}$ , the ideals $I,J$ are not nilpotent, and $\bigcap_k I^k=(0)$. Assume that $J$ is principal and  the ring $\mathcal B=\oplus_{m, n} J^m \cap I^n$ is noetherian. Then there exists a positive integer $t$ such that $v(m+t)=v(m)+v(t)$ for all $m \geq t$. 
\end{Proposition} 

\begin{proof} In the ring  $\oplus_{n \geq 0}I^n$ consider the filtration $\{F_m\}$ with $F_m=\oplus_{n \geq 0} J^m \cap I^n$. Since $\mathcal B=\oplus_{m \geq 0} F_m$ is noetherian, there exists a positive integer $t$ such that $F_{m+t}=F_m F_t$ for all $m \geq t$. We will prove that this implies $v(m+t)=v(m)+v(t)$ for all $m \geq t$. First note that the inequality $v(m+t) \geq v(m)+v(t)$ always holds. By contradiction, assume  that $v(m+t) >  v(m)+v(t)$ for some $m \geq t$. This implies that the component of degree $v(m)+v(t)+1$ in $F_{m+t}$ is  $J^{m+t}$, and since $F_{m+t}=F_m F_t$ we then obtain
$$ J^{m+t}=J^{t}(J^m \cap I^{v(m)+1}) + J^m (J^t \cap I^{v(t)+1}).$$  
Let $J=(z)$. Then we have
$$(z)^{m+t}=z^{m+t}(I^{v(m)+1}:z^m) + z^{m+t} (I^{v(t)+1}:z^t).$$
From the definition of $v(-)$, both $(I^{v(m)+1}:z^m)$ and $(I^{v(t)+1}:z^t)$ are contained in the maximal ideal, and  by the Nakayama Lemma, we must 
have $z$ nilpotent, contradicting our assumptions.
\end{proof}

\begin{Remark} 
\label{fields}
It is not always true that the ring $\mathcal B$ is noetherian. For such an example see \cite[Lemma 5.6]{Fields}. 
\end{Remark}

\medskip

Note that there are a few other natural conditions that ensure the
periodic increase  of the sequence $\{ v(m) \}_{m}$.  We comment on these below.

\begin{Remark}
\label{assocgraded}If the ring $\mathcal{G}(I)=\oplus_{n \geq 0} I^n/I^{n+1}$ is reduced, then we have $v(m)=mv(1)$ for all $m$. In particular, 
the sequence $v(m+1)-v(m)$ is constant.  Indeed, let $x \in J \setminus I^{v(1)+1}$. The image of $x$ in $I^{v(1)}/I^{v(1)+1} \subseteq \mathcal{G}(I)$ 
is nonzero, and since $\mathcal{G}(I)$ is reduced, so is the image of $x^m$ in $I^{mv(1)}/I^{mv(1)+1}$. This implies that $J^m \nsubseteq I^{mv(1)+1}$, and
hence $v(m) \leq m v(1)$.   
\end{Remark}

The point of view formulated in the above remark can be refined to include the case when
$J$ is not necessarily principal, but it comes at the expense of strengthening
the hypotheses.

\begin{Remark} Assume that $I$ is normal and $J=(a_1,\ldots,
  a_s)$. Then for every $m$ we have $v_I(J,m)=\min \{v_I((a_j),m)
  \mid j=1,\ldots,s \}$. Indeed, if $n:=\min \{v_I((a_j),m) \mid
  j=1,\ldots,s \}$, then $a_{j}^{m} \in I^n$ for all $j=1,\ldots,
  s$. This implies that $J^m \subseteq
  \overline{J^m}=\overline{(a_1^m,\ldots, a_s^m)} \subseteq
  \overline{I^n}=I^n$, so $v_I(J,m) \geq n$. On the other hand, if
  $v_I(J,m) >  n$, we have $J^m \subseteq I^{v_I((a_j),m)+1}$ for some
  $j$ and hence $a_j^m \in I^{v_I((a_j),m)+1}$, a contradiction. 
If $I$ is normal and all the rings $\oplus_{m, n} (a_j^m) \cap I^n$ are noetherian ($j=1,\ldots,s$), by Proposition~\ref{prop-per} we obtain that there exists $t_j$ such that $v_I((a_j), m+t_j)=v_I((a_j), m)+v_I((a_j),t_j)$ for $m \geq t_j$. If we have  $t_1=t_2=\ldots=t_s=t$ (the sequences $v_I((a_j), m)$ increase eventually in a periodic way with the same period), then we have $v_I(J, m+t)=v_I(J, m)+v_I(J,t)$ for $m \geq t$. Indeed, by the above observation, $v_I(J, m+t)= v_I((a_j), m+t_j)$ for some $j$, and hence $v_I(J, m+t)=v_I((a_j), m)+ v_I((a_j), t) \leq  v_I(J, m) +  v_I(J, t)$. The other inequality always holds.

Note that in the situation described in Remark~\ref{assocgraded}, when the associated graded ring  $\mathcal{G}(I)=\oplus_{n \geq 0} I^n/I^{n+1}$ is reduced (which implies that  $I$ is normal), we have  $t_1=t_2=\ldots=t_s=1.$
\end{Remark}

Our final observation introduces a bigraded ring associated to the
ideals $J$ and $I$ that can be used in examining the periodicity of the rate of change of  the sequence $\{ v(m) \}_m$.

\begin{Remark} Let $\mathcal C$ be the ring $\oplus _{m \geq 0, n\geq 0} F_{m,n}$, with $F_{m,n} = J^m \cap I^n/J^m \cap I^{n+1}$ and multiplication defined naturally such that $F_{m,n} F_{m',n'} \subseteq F_{m+m',n+n'}$. Let $F_m = \oplus_{n \geq 0} F_{m,n}$. Note that $F_m$ is a 
filtration on $\mathcal{G}(I)=\oplus_{n \geq 0} I^n/I^{n+1}$ and $F_{m,n} =0$ for $n< v(m)$, while $F_{m,v(m)} \neq 0$ for all $m$. As in the above remark, one can check that
$v(m+t)= v(m) +v(t)$ is equivalent to $F_{m,v(m)} F_{t,v(t)} \neq 0$.

So, if there exists $t$ such that $F_{t,v(t)}$ contains  a nonzerodivisor on $\mathcal C$, then $v(m+t)= v(m) +v(t)$ for all $m$. However, note that $\mathcal C$ a domain
implies that $F_0 = \mathcal{G}(I)$, the associated graded ring of $I$, is a domain as well, and then Remark~\ref{assocgraded} applies.
\end{Remark}

\section{Computations}

 In this section we describe a  method of determining $L_J(I)=\inf \{m/n \mid J^m \subseteq I^n \}$ (and  $l_I(J)=1/L_J(I)$) for two monomial ideals $I$ and $J$ in a polynomial ring $k[x_1,\ldots, x_r]$ over a field $k$. Whenever  $J=(a_1,\ldots, a_s)$, one has $L_J(I)=\max \{L_{(a_j)}(I) \mid j=1,\ldots, s \}$ ( \cite[Theorem 2]{Samuel1}), so we may assume that $J$ is a principal ideal. Let  $I=(x_1^{b_{i1}} x_2^{b_{i2}}\ldots x_r^{b_{ir}} \mid  i=1,\ldots,t)$ and $J=(x_1^{c_1}x_2^{c_2}\ldots x_r^{c_r})$.

First observe that  $J^m \subseteq I^n$ if and only if there exist nonnegative integers $y_1,\ldots,y_t$ with $y_1+\ldots+y_t=n$ such that
\begin{equation}
\sum_{i=1}^t b_{ij}y_i \leq c_j m \quad \text{for all} \quad j=1,\ldots,r.
\end{equation}
Set $B_{ij}=(1/c_j) b_{ij}$, $z_i=y_i/(y_1+\ldots+y_t)=y_i/n$ and $z=(z_1,\ldots,z_t) \in \mathbb{Q}^t$. 

So $ J^m \subseteq I^n$ if and only if there exist $z_i=y_i/n$ with $y_1+\ldots+y_t=n$ such that
\begin{equation}\label{eq4}
m/n \geq \frac{1}{nc_j} \sum_{i=1}^t b_{ij}y_i=\sum_{i=1}^t B_{ij}z_i \text{ for all }j=1, \ldots,r.
\end{equation}

Consider the function $\alpha: \mathbb{R}^t \to \mathbb{R}$ , $\alpha(z)=\max_{1 \leq j \leq r} \{ \sum_{i=1}^t B_{ij}z_i\}$  and  the subsets of the rationals $\Lambda_1=\{m/n \mid J^m \subseteq I^n \}$ and $\Lambda_2=\{\alpha(z) \mid  z_1, \ldots, z_t \in \mathbb{Q}_{\geq 0}, z_1+\ldots + z_t=1 \} $. We will prove that 
\begin{equation}\label{eq5}
\inf \Lambda_1 = \inf \Lambda_2
\end{equation}

The inequality $ \geq$ follows from (\ref{eq4}). For the other inequality, we will show that $\Lambda_2 \subseteq \Lambda_1$. Let $\alpha(z) \in \Lambda_2$ with $z_i=p_i/q$ ($1 \leq i \leq t$, $p_1+\ldots +p_t=q$, and $p_i, q$ nonnegative integers). The coefficients $B_{ij}$ are rationals, so  after clearing the denominators we obtain $\alpha(z)=h/lq$ for some nonnegative integers $h,l$.  By (\ref{eq4}), since $z_i=lp_i/lq$ for all $i$, we have $h/lq \in \Lambda_1$, which finishes the proof of (\ref{eq5}).

Note that $\inf \Lambda_2=\inf \{ \alpha(z) \mid  z_1, \ldots, z_t \in \mathbb{R}_{\geq 0}, z_1+\ldots + z_t=1 \}$, so we need  to minimize the function
$$\alpha(z)=\max \{ \sum_{i=1}^t B_{ij}z_i | j=1,\ldots,r \}$$ subject to the constraints
$$z_1,\ldots,z_t \geq 0 \quad \text{and} \quad z_1+\ldots+z_t=1.$$

Let $\displaystyle\Delta_k= \{z \in \mathbb R_{\geq 0}^t | \sum_{i=1}^t
B_{ik}z_i \geq \sum_{i=1}^t B_{ij}z_i  \text{ for all }  j \neq k \}
$. Clearly $\Delta_1 \cup \ldots \cup \Delta_r= \mathbb R_{\geq 0}^t$, so it is enough to minimize the function $\alpha$ on each $\Delta_k$.

In conclusion, for each $k=1,\ldots,r$,  the problem reduces to minimizing the objective function
$$\alpha(z)=\sum_{i=1}^t B_{ik}z_i$$ subject to the constraints
$$ z_1,\ldots, z_t \geq0, \quad z_1+\ldots+z_t=1 \quad \text{and}$$
$$ \sum_{i=1}^t B_{ik}z_i \geq \sum_{i=1}^t B_{ij}z_i \quad \text{for all} \quad j \neq k.$$ 

This is a classical problem linear programming problem which can be algorithmically solved using the simplex method.

\begin{Remark} In general, the limits $l_I(J)$ and $L_j(I)$ need not be reached by an element of the sequences  $\{v_I(J,m)/m\}_m$ and $\{w_J(I,n)/n\}_n$, respectively. However, in the monomial case, as the procedure described above shows, there exists a pair $(m,n)$ with $J^m \subseteq I^n$ and $L_J(I)=n/m$.

\end{Remark}

\begin{Example} Let $A=k[x,y]$ and $I=(x^3,x^2y, y^2)$, $J=(x^3y^7).$
In this case, $b_{11}=3, b_{12}=0, b_{21}=2, b_{22}=1, b_{31}=0, b_{32}=2$, $c_1=3, c_2=7$ and $B_{11}=3/3=1, B_{12}=0/7=0, B_{21}=2/3, B_{22}=1/7, B_{31}=0, B_{32}=2/7$. Then 
$$\Delta_1=\{(z_1,z_2,z_3) \in \mathbb R_{\geq 0}^3 \mid z_1+ (2/3) z_2 \geq (1/7) z_2 + (2/7) z_3 \}$$
and 
$$\Delta_2=\{(z_1,z_2,z_3) \in \mathbb R_{\geq 0}^3 \mid   (1/7) z_2 + (2/7) z_3 \geq z_1+
(2/3) z_2 \}.$$
By using a computer algebra system that has the simplex method
implemented, one can obtain that the minimum on each of the sets
$\Delta_1$ and $\Delta_2$ is $2/9$, and hence  $L_J(I)=2/9$.
\end{Example}
In fact, the minimum can occur only at the intersection of various regions $\Delta_k$ (in our case on $\Delta_1 \cap \Delta_2$),  for there are no critical points  in the interior of $\Delta_k$.

\begin{ack}The authors would like to thank Robert Lazarsfeld for a
  talk which inspired them to consider the problem treated in the
  article. They also thank Mel Hochster for pointing
  out to them the example mentioned in Remark~\ref{fields}.
\end{ack}

\end{document}